\documentclass[12pt]{article}

\usepackage{graphicx}
\usepackage{amsmath}
\usepackage{epsf}
\usepackage{graphpap}

\def\be{\begin{equation}}
\def\ee{\end{equation}}

\begin{document}
% \footnote{}

\begin{center}{\bf  \Large Comparison of Closed-form Solutions for the Lucas-Uzawa model via the Partial Hamitonian Approach and the Classical Approach}\\[2ex]
{R. Naz$^{a,*}$, Azam Chaudhry $^{b}$}\\[2ex]
{$^a$ Centre for Mathematics and Statistical Sciences,
  Lahore School of Economics, Lahore, 53200, Pakistan\\
  $^{b}$ Department of Economics,
Lahore School of Economics, Lahore, 53200, Pakistan\\
$^*$  Corresponding Author Email: drrehana@lahoreschool.edu.pk.\\
   Tel:
 92-3315439545\\
}
%
%\address{\\}
\end{center}
  \begin{abstract}
 In this paper we derive the closed-form solutions for the
Lucas-Uzawa growth model with the aid of the partial Hamiltonian
approach and then compare our results with those derived by the
classical approach \cite{chil}. The partial Hamiltonian approach
provides two first integrals \cite{naz2016} in the case where there
are no parameter restrictions and these two first integrals are
utilized to construct three sets of closed form solutions for all
the variables in the model.  First two first integrals are used to
find two closed form solutions, one of which is new to the
literature. We then use only one of the first integrals to derive a
third solution that is the same as that found in the previous
literature. We also show that all three solutions converge to the
same long run balance growth path.

\end{abstract}


\section{Introduction}
One of the foundations of modern economic growth theory is the
two-sector endogenous growth model developed by Lucas and Uzawa. The
model addresses the relationship between human capital accumulation
and economic growth and the idea behind the Lucas-Uzawa model (
\cite{lu} and \cite{uz}) is to determine optimal time paths for
consumption and the amount of labor devoted to the production of
capital in an economy which has constrained levels of physical and
human capital. One of the interesting features of the Lucas-Uzawa
model is the existence of multiple equilibria.

  The partial Hamiltonian approach \cite{naz} uses tools from
Lie group theory and is used to construct closed-form solutions of
dynamical systems such as those arising in economic growth theory.
This approach is unique and a significant departure from the rest of
the literature because unlike the previously used methods, the
partial Hamiltonian is applicable to an arbitrary system of ordinary
differential equations which means that it can be applied to more
complex models \cite{naz3}. In the context of our paper, the partial
Hamiltonian methodology yields a series of first integrals for a
system of ordinary differential equations and we use these first
integrals to find closed-form solutions for the Lucas-Uzawa model.
  In this paper we establish the closed-form solutions for the
Lucas-Uzawa model with the aid of the partial Hamiltonian approach
and we compare our results with those derived by the classical
approach \cite{chil}. The partial Hamiltonian approach provides two
first integrals \cite{naz2016} for case where there are no parameter
restrictions. We utilize these two first integrals to construct
closed form solutions for all variables of model for two different
scenarios: (i) $z=z^*$ and (ii) $z\not=z^*$ where
$z(t)=\frac{h(t)u(t)}{k(t)}$.

We begin by using both first integrals to construct closed-form
solutions for both scenarios. One solution is exactly the same as
derived by Chilarescu \cite{chil} and the second is completely new
to the literature. We then use only one first integral to determine
a different solution, again under fairly general parameter values,
and show that this is the same solution that has been derived in the
previous literature (see Chilarescu \cite{chil}).

We find that in the case where $z=z^*$, both the partial Hamiltonian
approach and the classical approach provide one solution. For the
$z\not=z^*$ case, the classical approach yields one solution while
the partial Hamiltonian approach yields the same solution as well as
providing one additional solution which is completely new to the
literature. What is especially interesting about our new solution is
that while the equilibrium levels of consumption and capital stock
in the new solution are equal to those found in the old solution,
the amount of labor allocated to the production of physical capital
and the level of human capital are different in the new solution.
The existence of three closed form solutions is new to the
literature and we also show that these closed form solutions all
converge to the same long run balanced growth path.

    It is important to mention here that under a specific parameter restriction $\sigma=\frac{\beta(\rho+\pi)}{2\pi\beta-\delta+\delta \beta-\pi}$ a
    third first integral was obtained and the closed form solution
for this case was new in the literature (see Naz et al
\cite{naz2016}).  The partial Hamiltonian approach provides three
solutions for the case in which there are no parameter restrictions.
For the case in which the specific parameter restriction
$\sigma=\frac{\beta(\rho+\pi)}{2\pi\beta-\delta+\delta \beta-\pi}$
is imposed, there is an additional solution as shown in Naz et al
\cite{naz}.  The classical approach provides only two closed-form
solutions for this model (see Chilarescu \cite{chil}) which means
that the partial Hamiltonian approach not only provides all the
solutions constructed in the previous literature (see Chilarescu
\cite{chil}) but also provides additional closed-form solutions
which are completely new to the literature.

The layout of the paper is as follows. In Section 2, we introduce
the Lucas-Uzawa model and provided expressions for the two first
integrals derived previously in \cite{naz}. In Section 3, the
closed-form solutions of the dynamical system of ODEs are
constructed by utilizing both first integrals in the case of no
parameter restrictions. In Section 4, we use only one first integral
to derive the closed form solutions for all the variables in the
model. A comparison of our results with those derived by the
classical approach is presented in Section 5. Finally, our
conclusions are summarized in Section 6.

\section{The Lucas-Uzawa model}

 The representative agent's utility
function is defined as
  \be Max_{c,u} \quad \int_0^{\infty}\frac{c^{1-\sigma}-1}{1-\sigma} e^{-\rho t} ,
  \; \sigma \not=1\label{(rn1)}\ee subject to the constraints
of  physical capital and human capital (see details of parameters:
\begin{eqnarray} \dot k(t)
= \gamma k^\beta u^{1-\beta} h^{1-\beta}-\pi k-c, \; k_0=k(0) \nonumber\\
\dot h(t) =\delta(1-u)h,\; h_0=h(0) \label{(rn2)}.\end{eqnarray}

where $1/{\sigma}$ is the constant elasticity of intertemporal
substitution, $\rho>0$ is the discount factor, $\beta$ is the
elasticity of output with respect to physical capital, $\gamma>0$ is
the technological levels in the goods sector, $\delta>0$ is the
technological levels in the education sector, $k$ is the level of
physical capital, $h$ is the level of human capital, $c$ is per
capita consumption and $u$ is the fraction of labor allocated to the
production of physical capital.

 The current value Hamiltonian function is defined as \be
H(t,c,k,\lambda)=\frac{c^{1-\sigma}-1}{1-\sigma}+\lambda[\gamma
k^\beta u^{1-\beta} h^{1-\beta} -\pi k-c]+\mu\delta(1-u)h,
\label{(rn3)}\ee where $\lambda(t)$ and $\mu(t)$ are costate
variables.  The transversality conditions are \be \lim_{t\to\infty}
e^{-\rho t}\lambda(t)k(t)=0 ,\; \lim_{t\to\infty} e^{-\rho
t}\mu(t)h(t)=0
 \label{(trans)}. \ee
 The Pontrygin's maximum principle provides following set of first order
 conditions:
 \be
\lambda=c^{-\sigma}, \label{(rn4)} \ee \be
u^{\beta}=\frac{\gamma(1-\beta)k^{\beta}h^{-\beta}}{\delta}\frac{\lambda}{\mu},\label{(rn5)}
\ee
 \be
\dot k(t) =\gamma k^\beta u^{1-\beta} h^{1-\beta}-\pi k-c,
\label{(rn6)} \ee
 \be \dot h(t) =\delta(1-u)h, \label{(rn7)} \ee
 \be   \dot
\lambda= -\lambda \gamma \beta
u^{1-\beta}k^{\beta-1}h^{1-\beta}+\lambda (\rho+\pi),
\label{(rn8)}\ee
 \be   \dot
\mu= \mu(\rho-\delta) .\label{(rn9)}\ee

The growth rates of consumption $c$ and physical capital $u$ are
given by \be \frac{\dot c}{c}= \frac{ \beta \gamma}{\sigma}
u^{1-\beta}k^{\beta-1}h^{1-\beta}-\frac{ \rho+\pi}{\sigma},
\label{(rn10)} \ee \be \frac{\dot u}{u}=
\frac{(\delta+\pi)(1-\beta)}{\beta}-\frac{ c}{k}+  \delta u.
\label{(rn11)} \ee
   The first integrals with no restriction on
parameters of economy derived via partial Hamiltonian approach by
Naz et al \cite{naz2016} are given by
\begin{eqnarray}
I_1=\frac{\gamma(1-\beta)}{\delta}c^{-\sigma}
{k}^{\beta}u^{-\beta}h^{-\beta }
e^{-(\rho-\delta )t},\nonumber\\
 I_2=\frac{c^{-\sigma}e^{-\rho t}}{1-\sigma}\bigg[ (\rho+\pi-\pi
\sigma)k-\sigma c-\beta \gamma (1-\sigma) (\frac{u h}{k})^{1-\beta}k
 \nonumber\\
 +\frac{(1-\beta)\gamma}{\delta }(\rho-\delta+\delta  \sigma)(\frac{k}{u
 h})^{\beta}h\bigg]. \label{(rn16)}
\end{eqnarray}

Another first integral exists under parameter restriction
$\sigma=\frac{\beta(\rho+\pi)}{2\pi\beta-\delta+\delta \beta-\pi}$
provided $2\pi\beta-\delta+\delta \beta-\pi>0$ to ensure that
$\sigma>0$. The complete analysis in terms of closed-from solution
under parameter restriction
$\sigma=\frac{\beta(\rho+\pi)}{2\pi\beta-\delta+\delta \beta-\pi}$
was provided by Naz et al \cite{naz2016}.

\section{Closed-form
solution for the Lucas-Uzawa model under fairly general values of
parameters via $I_1$ and $I_2$}

The closed-form solution via $I_1$ and $I_2$ for the case
$z\not=z^*$ case was provided in Naz et al \cite{naz2016} and we
provide this solution in simplified form. Here, we provide closed
form solutions for the original variables $c(t), u(t), k(t), h(t),
\lambda(t)$ and $\mu(t)$ for the case $z\not=z^*$ explicitly in
terms of variable $z(t)$. Moreover, we provide closed form solutions
for the original variables $c(t), u(t), k(t), h(t), \lambda(t)$ and
$\mu(t)$ for the case $z=z^*$ as well.

 By setting
$I_1=c_1$, we obtain \be \frac{\gamma(1-\beta)}{\delta}c^{-\sigma}
{k}^{\beta}u^{-\beta}h^{-\beta } e^{-(\rho-\delta )t}=c_1,
\label{(f1)} \ee where $c_1$ is an arbitrary constant. Introducing
$z=\frac{h u}{k}$, Equation {(\ref {(f1)})}  can be re-written as
\be
z=\bigg(\frac{(1-\beta)\gamma}{c_1\delta}\bigg)^{\frac{1}{\beta}}\lambda^{\frac{1}{{\beta}}}
e^{-\frac{(\rho-\delta)}{\beta}t}, \label{(f2)} \ee where
$\lambda=c^{-\sigma}$.
  Equation {(\ref {(rn8)})} with the aid of
Equation {(\ref {(f2)})} yields following Bernoulli's differential
equation for $\lambda$
 \be   \dot
\lambda-\lambda (\rho+\pi)= - \beta \gamma
\bigg[\frac{(1-\beta)\gamma}{c_1\delta}\bigg]^{\frac{1-\beta}{\beta}}\lambda^{\frac{1}{{\beta}}}
e^{-\frac{(\rho-\delta )(1-\beta)}{\beta}t}, \label{(f3)}\ee and
thus we have \be \lambda=c^{-\sigma}= \bigg[\frac{ \beta
\gamma^{\frac{1}{\beta}}}{\delta +\pi}
\bigg(\frac{1-\beta}{c_1\delta}\bigg)^{\frac{1-\beta}{\beta}}
e^{-\frac{(\rho-\delta )(1-\beta)}{\beta}t} +c_2
e^{-\frac{(1-\beta)(\rho+\pi)}{\beta}t}
\bigg]^{\frac{\beta}{\beta-1}},\label{(f4)}\ee where $c_2$ is an
arbitrary constant.  We found explicit solution for $z$ after
substituting value of $\lambda$ from Equation {(\ref {(f4)})} into
Equation {(\ref {(f2)})}
 \be
z= \bigg[\frac{\beta \gamma}{\delta +\pi}
 +c_2
\bigg(\frac{(1-\beta)\gamma}{c_1\delta}\bigg)^{\frac{\beta-1}{\beta}}e^{-\frac{(1-\beta)(\delta
+\pi)}{\beta}t} \bigg]^{\frac{1}{\beta-1}},\label{(oldf5)}\ee and $
z^*=\bigg(\frac{\beta \gamma}{\delta +\pi} \bigg
)^{\frac{1}{\beta-1}}$ is the steady state solution. It is worthy to
mention here that the system of differential equations {(\ref
{(rn6)})}-{(\ref {(rn9)})} provides two sets of solutions depending
on $c_2=0$ and $c_2\not=0$. Thus  we discuss two scenarios \\
  Scenario I:  $c_2=0$  and thus $z=z^*$\\
   Scenario II: $c_2\not= 0$  and thus $z\not=z^*$.
\subsection{Scenario I:  $c_2=0$  and thus $z=z^*$}
  For
  $c_2=0$, we have

\be z(t)= \bigg(\frac {\delta +\pi}{\beta \gamma}
  \bigg)^{\frac{1}{1-\beta}}=z^*.\label{(f500)}\ee
  Equation {(\ref {(f2)})} yields
    \be \lambda(t)=\frac{c_1
\delta}{(1-\beta)\gamma}e^{(\rho-\delta )t}z^{
* \beta},\label{(f501)}\ee
and thus \be c(t)=\bigg(\frac{(1-\beta)\gamma z^{* \beta}}{c_1
\delta }\bigg)^{\frac{1}{\sigma}}e^{-\frac{(\rho-\delta
)}{\sigma}t}. \label{(f502)}\ee Using initial condition $c(0)=c_0$
and $k(0)=k_0$, Equations {(\ref {(f502)})} gives
$$c_0=\bigg(\frac{(1-\beta)\gamma }{c_1 \delta z^{* \beta}
}\bigg)^{\frac{1}{\sigma}}.$$

 Equation {(\ref {(rn6)})} provides
following solution for $k(t)$: \be
 k(t)=\frac{\sigma \beta
c_0}{ \sigma \pi+\delta (\sigma-\beta)-(\pi
\sigma-\rho)\beta}e^{-\frac{(\rho-\delta )}{\sigma}t}+a_1
e^{-\frac{(\pi \beta-\pi-\delta )}{\beta}t}. \label{(ff502)}\ee The
transversality condition {(\ref {(trans)})} for $k$ is satisfied
provided  $\delta <\rho+\delta \sigma$ and $a_1=0$. Using initial
condition $k(0)=k_0$, we have \be
 \frac{c_0}{k_0}=\frac{\delta+\pi(1-\beta)}{\beta}-\frac{\delta-\rho}{\sigma}>0, \label{(ff604)}\ee
as $\frac{\delta+\pi-\pi \beta}{\beta}-\frac{\delta-\rho}{\sigma}
>0$ (see proof in Proposition 1 \cite{chil}). Equations {(\ref
{(rn5)})} and
 {(\ref {(f1)})} provide following expression of costate
variable $\mu$ \be \mu= c_1 e^{(\rho-\delta )t}.\label{(f499)}\ee

 Next, we set $I_2=c_3$ and after
some simplifications, we arrive at following expression for $h$
\begin{eqnarray}
h(t)=\frac{\delta}{\gamma(1-\beta)(\rho-\delta+\delta
\sigma)}z^{*\beta}\bigg[c_3 (1-\sigma)c_0^{\sigma}e^{-\delta
t}\nonumber\\+\bigg(\beta
\gamma(1-\sigma)z^{*1-\beta}k_0-(\rho+\pi-\pi\sigma)k_0 +\sigma
c_0\bigg)e^{-\frac{(\rho-\delta)}{\sigma}t}\bigg]\end{eqnarray}

The transversality condition {(\ref {(trans)})} for $h$ is satisfied
provided $c_3=0$ and  $\delta <\rho+\delta \sigma$. The initial
condition $h(0)=h_0$ yields \be
h(t)=h_0e^{-\frac{(\rho-\delta)}{\sigma}t} \ee where \be
h_0=\frac{z^* k_0 \delta \sigma}{\rho-\delta+\delta \sigma}.\ee
Finally, $u=\frac{z k}{h}$ gives \be u=\frac{\rho-\delta+\delta
\sigma}{\delta \sigma},\ee
 and this completes the solution.
 We can summarize these solutions for all variables in the following simple
forms: \begin{eqnarray}  c(t)=c_0e^{-\frac{(\rho-\delta
)}{\sigma}t}, \nonumber\\
 k(t)=k_0e^{-\frac{(\rho-\delta )}{\sigma}t},\nonumber\\
  u(t)=\frac{\rho-\delta+\delta
\sigma}{\delta \sigma}=u^*,\label{(sol1)} \\
 h(t)=h_0 e^{-\frac{(\rho-\delta )}{\sigma
}t}, \nonumber\\
\lambda(t)=c_0^{-\sigma}e^{(\rho-\delta)t},\nonumber\\
 \mu= c_1  e^{(\rho-\delta)
 t},\nonumber\\
 z^*=\bigg(\frac {\delta +\pi}{\beta \gamma}
  \bigg)^{\frac{1}{1-\beta}} \nonumber,\end{eqnarray}
  provided $ \delta <\rho+\delta \sigma ,
c_0=\bigg(\frac{(1-\beta)\gamma }{c_1 \delta z^{* \beta}
}\bigg)^{\frac{1}{\sigma}},
 \frac{c_0}{k_0}=\frac{\delta+\pi(1-\beta)}{\beta}-\frac{\delta-\rho}{\sigma}>0,h_0=\frac{z^* k_0}{u^*}$.

\subsection{Scenario II: $c_2\not= 0$  and thus $z\not=z^*$}
The expression for $\lambda$ given in Equation {(\ref {(f4)})} can
be alternatively given as \be \lambda=\frac{c_1
\delta}{(1-\beta)\gamma}e^{(\rho-\delta)t}z^{\beta}.
\label{(ff2)}\ee Euqtaion {(\ref {(oldf5)})} for $z(t)$, with
initial condition
 $z(0)=z_0$ takes following form:

 \be z(t)= \frac{z^*z_0}{[(z^{*1-\beta}-z_0^{1-\beta})e^{-\frac{(1-\beta)(\delta
+\pi)}{\beta}t}+z_0^{1-\beta}]^{\frac{1}{1-\beta}} },\label{(f5)}\ee
where $c_2
=\frac{z_0^{1-\beta}-z^{*1-\beta}}{(\frac{(1-\beta)\gamma}{c_1\delta})^{\frac{\beta-1}{\beta}}}$.

The variable $c(t)$ with $c(0)=c_0$ takes following form:\be c=c_0
z_0^{\frac{\beta}{\sigma}}e^{-\frac{(\rho-\delta
)}{\sigma}t}z^{-\frac{\beta}{\sigma}},\;c_0
z_0^{\frac{\beta}{\sigma}}=\bigg(\frac{c_1
\delta}{(1-\beta)\gamma}\bigg)^{-\frac{1}{\sigma}} \label{(fff2)}\ee
where $z$ is same as given in {(\ref {(f5)})}. The differential
equation {(\ref {(rn6)})} for $k$ results in following integrable
differential equation \be \dot k+(\pi -\gamma z^{1-\beta}) k=-c_0
z_0^{\frac{\beta}{\sigma}}e^{-\frac{(\rho-\delta
)}{\sigma}t}z^{-\frac{\beta}{\sigma}},\label{(f6)} \ee  and it
provides
 \be
k(t)=\bigg(a_3-\frac{c_0 z_0^{\frac{\beta}{\sigma}}}{(\pi+\delta
)^{\frac{1}{1-\beta}}}F(t) \bigg)(\pi+\delta
)^{\frac{1}{1-\beta}}z(t)^{-1}e^{\frac{(\pi+\delta -\pi
\beta)}{\beta}t},\label{(950)}\ee where \be F(t)=\int_0^t
z(t)^{\frac{\sigma-\beta}{\sigma}} e^{-(\frac{\delta+\pi-\pi
\beta}{\beta}-\frac{\delta-\rho}{\sigma} )t} dt\ee $a_3$ is
arbitrary constant of integration and  $\frac{\delta+\pi-\pi
\beta}{\beta}-\frac{\delta-\rho}{\sigma} >0$ (see proof in
Proposition 1 \cite{chil}). The initial condition $k(0)=k_0$ yields
$a_3=\frac{k_0z_0}{(\pi+\delta)^{\frac{1}{1-\beta}}}$ and thus
expression for $k(t)$ simplifies to the following form: \be
k(t)=\bigg(\frac{k_0}{c_0z_0^{\frac{\beta-\sigma}{\sigma}}}-F(t)
\bigg) c_0z_0^{\frac{\beta}{\sigma}}z(t)^{-1}e^{\frac{(\delta+\pi
-\pi \beta)}{\beta}t}.\label{(950a)}\ee
  The
transversality condition {(\ref {(trans)})} for $k$ is satisfied
provided  $\delta <\rho+\delta \sigma$ and
 \be \lim_{t\to\infty} F(t)=\frac{k_0}{c_0z_0^{\frac{\beta-\sigma}{\sigma}}}.\label{(950aa)}\ee
It is important to mention here that the integrand of $F(t)$ is
positive and bounded therefore $\lim_{t\to\infty} F(t)$ is a finite
number.  Equation {(\ref {(rn5)})} provides following expression for
the costate variable \be \mu= c_1 e^{(\rho-\delta )t}.\label{(m)}\ee

 Setting $I_2=c_3$, we find
\begin{eqnarray}\frac{c^{-\sigma}e^{-\rho t}}{1-\sigma}\bigg[
(\rho+\pi-\pi \sigma)k-\sigma c-\beta \gamma (1-\sigma) (\frac{u
h}{k})^{1-\beta}k \nonumber\\
 +\frac{(1-\beta)\gamma}{\delta }(\rho-\delta+\delta  \sigma)(\frac{k}{u
 h})^{\beta}h\bigg]=c_3 \label{(f7)} \end{eqnarray}
 and this gives
 \be h(t)= \frac{\left( {\frac {c_{{3}} \left( 1-\sigma \right) }{\lambda\,{{\rm e}^{-
\rho\,t}}}}- \left( \rho+\pi -\pi \,\sigma \right) k+\sigma
c+\beta\,\gamma\, \left( 1-\sigma \right) {z}^{1-\beta }k \right)
\delta  {z}^{\beta}}{ \gamma(1-\beta)(\rho-\delta  +\delta  \sigma)
}.
 \label{(f8)} \ee
 The transversality condition for $h$ requires to choose $c_3=0$ and
  $\delta <\rho+\delta \sigma $.  Finally, the variable $u$ can be
determined from $u={z k}/{h}$ and it simplifies to \begin{eqnarray}
u(t)=\frac{\delta^{-1}\gamma(1-\beta)(\rho-\delta+\delta
\sigma)(\frac{k_0}{c_0
z_0^{\frac{\beta-\sigma}{\sigma}}}-F(t))}{[\beta
\gamma(1-\sigma)-(\rho+\pi-\pi \sigma)z^{\beta-1}](\frac{k_0}{c_0
z_0^{\frac{\beta-\sigma}{\sigma}}}-F(t))+\sigma
z^{\beta-\frac{\beta}{\sigma}e^{-(\frac{\delta+\pi-\pi
\beta}{\beta}-\frac{\delta-\rho}{\sigma} )t}}}.\nonumber
\end{eqnarray}
 This completes
the solution. We then apply the initial conditions $h(0)=h_0$,
$u_0=\frac{z_0k_0}{h_0}$ and we summarize the solutions for all
variables as follows:
\begin{eqnarray} c(t)=c_0
z_0^{\frac{\beta}{\sigma}}e^{-\frac{(\rho-\delta
)}{\sigma}t}z^{-\frac{\beta}{\sigma}}, \nonumber\\
k(t)=\bigg(\frac{k_0}{c_0z_0^{\frac{\beta-\sigma}{\sigma}}}-F(t)
\bigg) c_0z_0^{\frac{\beta}{\sigma}}z(t)^{-1}e^{\frac{(\delta
+\pi-\pi \beta)}{\beta}t}
,\nonumber\\
 h(t)=\frac{h_0}{z_0[\sigma c_0 z_0^{\beta-1}
-(\rho+\pi-\pi \sigma)k_0 z_0^{\beta-1}+\beta \gamma(1-\sigma)k_0]}[
\sigma c_0 z_0^{\frac{\beta}{\sigma}}e^{-\frac{(\rho-\delta
)}{\sigma}t}z^{-\frac{\beta}{\sigma}+\beta} \nonumber\\
+(\beta \gamma (1-\sigma) -(\rho+\pi-\pi \sigma)z^{\beta-1})
(\frac{k_0}{c_0z_0^{\frac{\beta-\sigma}{\sigma}}}-F(t) )
c_0z_0^{\frac{\beta}{\sigma}}e^{\frac{(\delta+\pi
-\pi \beta)}{\beta}t}], \nonumber\\
 u(t)=\frac{u_0}{k_0}[\sigma c_0 z_0^{\beta-1}-(\rho+\pi-\pi \sigma)k_0 z_0^{\beta-1}+\beta \gamma(1-\sigma)k_0]\nonumber\\
 \times \frac{(\frac{k_0}{c_0
z_0^{\frac{\beta-\sigma}{\sigma}}}-F(t))}{[\beta
\gamma(1-\sigma)-(\rho+\pi-\pi \sigma)
 z^{\beta-1}](\frac{k_0}{c_0 z_0^{\frac{\beta-\sigma}{\sigma}}}-F(t))+\sigma z^{\beta-\frac{\beta}{\sigma}}e^{-(\frac{\delta+\pi-\pi
\beta}{\beta}-\frac{\delta-\rho}{\sigma} )t}},\nonumber\\
 \lambda(t)=\frac{c_1
\delta}{(1-\beta)\gamma}e^{(\rho-\delta)t}z^{\beta},\nonumber\\
 \mu(t)= c_1 e^{(\rho-\delta )t}
\nonumber,\end{eqnarray} where \begin{eqnarray} F(t)=\int_0^t
z(t)^{\frac{\sigma-\beta}{\sigma}} e^{-(\frac{\delta+\pi-\pi
\beta}{\beta}-\frac{\delta-\rho}{\sigma} )t} dt,\nonumber\\ z(t)=
\frac{z^*z_0}{[(z^{*1-\beta}-z_0^{1-\beta})e^{-\frac{(1-\beta)(\delta
+\pi)}{\beta}t}+z_0^{1-\beta}]^{\frac{1}{1-\beta}} },\label{(sol2)}\\
\lim_{t\to\infty}
F(t)=\frac{k_0}{c_0z_0^{\frac{\beta-\sigma}{\sigma}}},\nonumber\\
 \rho<\delta <\rho+\delta \sigma, \frac{\delta+\pi-\pi
\beta}{\beta}-\frac{\delta-\rho}{\sigma}
>0,\nonumber\\
c_0 z_0^{\frac{\beta}{\sigma}}=\bigg(\frac{c_1
\delta}{(1-\beta)\gamma}\bigg)^{-\frac{1}{\sigma}},\nonumber\\
\frac{\gamma(1-\beta)(\rho-\delta+\delta \sigma)}{\delta}\nonumber\\=\frac{u_0}{k_0}[\sigma c_0 z_0^{\beta-1}-(\rho+\pi-\pi \sigma)k_0 z_0^{\beta-1}+\beta \gamma(1-\sigma)k_0],\nonumber\\
 z^*=\bigg(\frac{\beta \gamma}{\delta +\pi} \bigg
)^{\frac{1}{\beta-1}}\nonumber.\end{eqnarray}
 The first integrals $I_1$ and $I_2$ yield two solutions of the dynamical
system of ODEs {(\ref {(rn6)})}-{(\ref {(rn11)})} given in equations
{(\ref {(sol1)})} and {(\ref {(sol2)})} for fairly general values of
$\sigma$ and $\beta$. It is straight forward to show that \be
\lim_{t\to\infty} u(t)=u^*.\ee
\section{Closed-form
solution for the Lucas-Uzawa model under fairly general values of
parameters via $I_1$} Now we show how one can utilize only one first
integral $I_1$ to derive the same closed-from solution as in the
existing literature which was also derived by Chilarescu \cite{chil}
via the classical approach. By setting $I_1=a_1$, we will arrive at
equations {(\ref {(f1)})}-{(\ref {(oldf5)})}. Here also the
following two scenarios will
arise:\\

 Case I:  $a_2=0$ \\
Case II: $a_2\not= 0$ .\\

Equations (\ref{(f500)})-(\ref{(f499)}) for variables $c(t), k(t),
\lambda(t)$ and $\mu(t)$ will follow from the previous section. Now
instead of utilizing the second first integral, we proceed as
follows to derive closed form solutions for $h(t)$ and $u(t)$.

Equation {(\ref {(rn11)})} simplifies to \be \frac{\dot u}{u}=
\frac{\delta-\rho-\delta \sigma}{\sigma}+ \delta u. \label{(rnnn11)}
\ee Equation {(\ref {(rnnn11)})} gives
  \be u(t)=
\frac{\frac{\delta-\rho-\delta
\sigma}{\sigma}}{a_2(\frac{\delta-\rho-\delta
\sigma}{\sigma})e^{-\frac{\delta-\rho-\delta \sigma}{\sigma}
t}-\delta }. \label{(rn701)} \ee  Finally, $h=\frac{z k}{u}$ gives
\be
 h(t)=
\frac{a_2(\frac{\delta-\rho-\delta
\sigma}{\sigma})e^{-\frac{\delta-\rho-\delta \sigma}{\sigma}
t}-\delta }{\frac{\delta-\rho-\delta
\sigma}{\sigma}}z^*k_0e^{-\frac{(\rho-\delta)}{\sigma}t}.
\label{(ff802)}\ee The transversality condition {(\ref {(trans)})}
for $h$ is satisfied provided $\delta <\rho+\delta \sigma$   and
$a_2=0$.  Thus we arrive at the same solution for all variables as
given by Equation {(\ref{(sol1)})} .

 For $a_2\not= 0$, we will follow same procedure as described in previous section to derive equations {(\ref {(ff2)})}-{(\ref
 {(m)})}. Substituting $c(t)$ and $k(t)$ from equations {(\ref
 {(fff2)})} and {(\ref {(950a)})} into equation {(\ref {(rn11)})}, we
 have
 \be \frac{\dot
u}{u}= \frac{(\delta+\pi)(1-\beta)}{\beta}-\frac{F'(t)}{ \frac{k_0
}{c_0 z_0^{\frac{\beta-\sigma}{\sigma}}}- F(t)}+\delta u.
\label{(rnn11)} \ee  The solution of equation {(\ref {(rnn11)})}
with initial condition $u(0)=u_0$ is given by

\be  u(t)=\frac{\frac{(\delta+\pi)(1-\beta)}{\beta} u_0[\frac{k_0
}{c_0 z_0^{\frac{\beta-\sigma}{\sigma}}}-
F(t)]}{[(\frac{(\delta+\pi)(1-\beta)}{\beta}  +\delta u_0)\frac{k_0
}{c_0 z_0^{\frac{\beta-\sigma}{\sigma}} }-\delta u_0 G(t)]e^{
-\frac{(\delta+\pi)(1-\beta)}{\beta}t}-\delta u_0 [\frac{k_0 }{c_0
z_0^{\frac{\beta-\sigma}{\sigma}}}- F(t)]}, \label{(u3)}\ee where
\be G(t)=\int_0^t z(t)^{\frac{\sigma-\beta}{\sigma}}
e^{-\frac{\delta \sigma-\delta-\rho}{\sigma} t} dt.\ee The solution
{(\ref {(u3)})} holds provided \be \lim_{t\to\infty}
\bigg[(\frac{(\delta+\pi)(1-\beta)}{\beta}  +\delta u_0)\frac{k_0
}{c_0 z_0^{\frac{\beta-\sigma}{\sigma}} }-\delta u_0 G(t) \bigg]=0,
\label{(chil3)} \ee with \be\lim_{t\to\infty} G(t)=
\frac{(\frac{(\delta+\pi)(1-\beta)}{\beta}  +\delta u_0)}{\delta
u_0}\lim_{t\to\infty} F(t),\ee  and $\lim_{t\to\infty} F(t)$ is
given in (\ref{(950aa)}). The variable $h$ can be determined from
$h={z k}/{u}$ and is given by

\begin{eqnarray}h(t)=\bigg[(\frac{(\delta+\pi)(1-\beta)}{\beta}  +\delta
u_0)\frac{k_0 }{c_0 z_0^{\frac{\beta-\sigma}{\sigma}} }-\delta u_0
G(t)]e^{ -\frac{(\delta+\pi)(1-\beta)}{\beta}t}\nonumber\\
-\delta u_0 [\frac{k_0 }{c_0 z_0^{\frac{\beta-\sigma}{\sigma}}}-
F(t)\bigg]\times
\frac{c_0z_0^{\frac{\beta}{\sigma}}}{\frac{(\delta+\pi)(1-\beta)}{\beta}
u_0}e^{\frac{(\pi+\delta -\pi \beta)}{\beta}t},
\label{(h3)}\end{eqnarray}

  We can summarize the solutions for all variables as follows:
\begin{eqnarray} c(t)=c_0
z_0^{\frac{\beta}{\sigma}}e^{-\frac{(\rho-\delta
)}{\sigma}t}z^{-\frac{\beta}{\sigma}}, \nonumber\\
k(t)=\bigg(\frac{k_0}{c_0z_0^{\frac{\beta-\sigma}{\sigma}}}-F(t)
\bigg) c_0z_0^{\frac{\beta}{\sigma}}z(t)^{-1}e^{\frac{(\pi+\delta
-\pi \beta)}{\beta}t}
,\nonumber\\
h(t)=[(\frac{(\delta+\pi)(1-\beta)}{\beta}  +\delta u_0)\frac{k_0
}{c_0 z_0^{\frac{\beta-\sigma}{\sigma}} }-\delta u_0
G(t)]e^{ -\frac{(\delta+\pi)(1-\beta)}{\beta}t}\nonumber\\
-\delta u_0 [\frac{k_0 }{c_0 z_0^{\frac{\beta-\sigma}{\sigma}}}-
F(t)]\times
\frac{c_0z_0^{\frac{\beta}{\sigma}}}{\frac{(\delta+\pi)(1-\beta)}{\beta}
u_0}e^{\frac{(\pi+\delta -\pi \beta)}{\beta}t},
\nonumber \\
u(t)=\frac{\frac{(\delta+\pi)(1-\beta)}{\beta} u_0[\frac{k_0 }{c_0
z_0^{\frac{\beta-\sigma}{\sigma}}}-
F(t)]}{[(\frac{(\delta+\pi)(1-\beta)}{\beta}  +\delta u_0)\frac{k_0
}{c_0 z_0^{\frac{\beta-\sigma}{\sigma}} }-\delta u_0 G(t)]e^{
-\frac{(\delta+\pi)(1-\beta)}{\beta}t}-\delta u_0 [\frac{k_0 }{c_0
z_0^{\frac{\beta-\sigma}{\sigma}}}- F(t)]}\nonumber\\
 \lambda(t)=\frac{c_1
\delta}{(1-\beta)\gamma}e^{(\rho-\delta)t}z^{\beta},\nonumber\\
 \mu(t)= c_1 e^{(\rho-\delta )t}
\nonumber,\end{eqnarray}  where

\begin{eqnarray} \rho<\delta <\rho+\delta \sigma, \frac{\delta+\pi-\pi
\beta}{\beta}-\frac{\delta-\rho}{\sigma}
>0,\nonumber\\ F(t)=\int_0^t
z(t)^{\frac{\sigma-\beta}{\sigma}} e^{-(\frac{\delta+\pi-\pi
\beta}{\beta}-\frac{\delta-\rho}{\sigma} )t} dt, \nonumber\\
G(t)=\int_0^t z(t)^{\frac{\sigma-\beta}{\sigma}} e^{-\frac{\delta
\sigma-\delta+\rho}{\sigma} t} dt,
\label{(sol3)} \\
z(t)=
\frac{z^*z_0}{[(z^{*1-\beta}-z_0^{1-\beta})e^{-\frac{(1-\beta)(\delta
+\pi)}{\beta}t}+z_0^{1-\beta}]^{\frac{1}{1-\beta}} },\nonumber\\
 c_0
z_0^{\frac{\beta}{\sigma}}=\bigg(\frac{c_1
\delta}{(1-\beta)\gamma}\bigg)^{-\frac{1}{\sigma}}, \nonumber\\
\lim_{t\to\infty}
F(t)=\frac{k_0}{c_0z_0^{\frac{\beta-\sigma}{\sigma}}},\nonumber\\
\lim_{t\to\infty} \bigg[(\frac{(\delta+\pi)(1-\beta)}{\beta} +\delta
u_0)\frac{k_0 }{c_0 z_0^{\frac{\beta-\sigma}{\sigma}} }-\delta u_0
G(t) \bigg]=0,\nonumber\\ \lim_{t\to\infty} G(t)=
\frac{(\frac{(\delta+\pi)(1-\beta)}{\beta}  +\delta u_0)}{\delta
u_0}\lim_{t\to\infty} F(t), \nonumber\\z^*=\bigg(\frac{\beta
\gamma}{\delta +\pi} \bigg
)^{\frac{1}{\beta-1}}.\nonumber\end{eqnarray} The closed form
solution {(\ref {(sol3)})} derived by only utilizing $I_1$ is
exactly the same as derived by Chilarescu \cite{chil} via the
classical approach.

\section{Comparison of Closed-form Solutions for the Lucas-Uzawa model}

The closed form solutions of the Lucas-Uzawa model have been derived
in the literature by using both the newly developed partial
Hamiltonian approach and the classical approach.  The partial
Hamiltonian approach utilizes Lie group theoretical techniques to
construct a closed-form solution. Using this partial Hamiltonian
methodology, we have established three sets of closed form solutions
{(\ref {(sol1)})}, {(\ref {(sol2)})} and {(\ref {(sol3)})} for the
Lucas-Uzawa model with no parameter restrictions. For the $z=z^*$
case only one solution arises which is given in {(\ref {(sol1)})}
whereas for the $z\not=z^*$ case we obtained two solutions which are
given in {(\ref {(sol2)})} and {(\ref {(sol3)})}.

We have shown how first integrals derived via the partial
Hamiltonian approach can be utilized to construct multiple closed
form solutions. It is shown that for the case where $z=z^*$, the
closed form solution {(\ref {(sol1)})} which was previously derived
in the literature using the classical approach (see Chilarescu
\cite{chil}) can also be constructed by utilizing either one first
integral $I_1$ or by using two first integrals $I_1$ and $I_2$. For
the case where $z\not=z^*$, we have obtained the second solution
{(\ref {(sol3)})} found by Chilarescu \cite{chil} and this solution
was obtained by utilizing only one first integral $I_1$. Using our
partial Hamiltonian methodology, we also arrive at an additional
solution {(\ref {(sol2)})} in the case where $z\not=z^*$ by
utilizing the two first integrals $I_1$ and $I_2$.

Thus there exists three sets of closed form solutions {(\ref
{(sol1)})}, {(\ref {(sol2)})} and {(\ref {(sol3)})} for the
Lucas-Uzawa model with no parameter restrictions.  Chilarescu
\cite{chil} utilized the classical approach to provide two solutions
with no restrictions on parameters which were the same as the
solutions in {(\ref {(sol1)})} and {(\ref {(sol3)})} obtained using
the partial Hamiltonian approach. The partial Hamiltonian approach
also provided one additional solution given in {(\ref {(sol2)})}
which is new to the literature.

In the long run all these solutions converge to the same steady
state which is important in the context of economic growth theory
and gives rise to multiple equilibria.  For the $z=z^*$ case, it is
straightforward to compute this simple closed form solution by any
of the existing techniques in the literature. At the same time, the
case in which $z\not=z^*$ is also important. A comparison of the
closed-form solutions {(\ref {(sol2)})} and {(\ref {(sol3)})} in the
second case shows that the expressions for consumption $c$, physical
capital stock $k$, and the costate variables $\lambda$ and $\mu$ are
the same in both solutions. On the other hand, the expressions for
the fraction of labor devoted to physical capital, $u$, and the
level of human capital, $h$, are different in our newly obtained
solution. Another important result is that the previously obtained
closed form solution {(\ref {(sol3)})} involves two numerically
computable functions $F(t)$ and $G(t)$ whereas our newly closed form
solution {(\ref {(sol2)})} involves only one numerically computable
function $F(t)$. So our newly derived closed form solution {(\ref
{(sol2)})} which was obtained from the partial Hamiltonian approach
is fundamentally different and is also in a simpler form than the
previously obtained solution {(\ref {(sol3)})}.

 It is straight forward to show, for the closed form solutions
{(\ref {(sol1)})}, {(\ref {(sol2)})} and {(\ref {(sol3)})}, that \be
\lim_{t\to\infty} u(t)=u^*.\ee It is worthy to mention here that
l'H$\hat{o}$pital rule is applied to establish $\lim_{t\to\infty}
u(t)$.
 The growth rates of consumption, $c$,
physical capital, $k$, and human capital, $h$, decrease over time
and approach $\frac{\delta-\rho}{\sigma}$  as $t\mapsto\infty$.
Also, the growth rate of the fraction of labor allocated to the
production of physical capital, $u$, approaches zero as
$t\mapsto\infty$ whereas the growth rates of both costate variables,
$\lambda$ and $\mu$, equal $(\rho-\delta)$ as $t\mapsto\infty$.

 {\bf Remark 1:} It
is important to mention here that partial Hamiltonian approach
provided a new specific restriction
$\sigma=\frac{\beta(\rho+\pi)}{2\pi\beta-\delta+\delta \beta-\pi}$
on parameters for which another first integral was established. Naz
et al \cite{naz2016} provided detailed analysis of closed-form
solution under this specific restriction. The closed form solution
under this specific restriction was not derived before in
literature.

 {\bf Remark 2:} It is worthy to mention here that one cannot
 derived closed-form solutions for the case $\sigma=\beta$ directly
 from {(\ref {(sol1)})}, {(\ref {(sol2)})} and {(\ref {(sol3)})}. We
 need to derive solutions by utilizing $I_1$ and $I_2$ or only
 $I_1$.

\section{Conclusions}
In this paper we establish closed-form solutions for the Lucas-Uzawa
model with the aid of the partial Hamiltonian approach and then
compare our results with those derived by the classical approach.
The partial Hamiltonian approach provides two first integrals
\cite{naz2016} in the case where there are no parameter
restrictions. We utilize these two first integrals to construct
closed form solutions for all the variables in the model for two
different scenarios: (i) $z=z^*$ and (ii) $z\not=z^*$ where
$z(t)=\frac{h(t)u(t)}{k(t)}$.

 For the $z=z^*$ case, both the
partial Hamiltonian approach and the classical approach provide one
solution. For the $z\not=z^*$ case, the classical approach yields
one solution while the the partial Hamiltonian approach provides the
same solution as well as one additional solution which is completely
new to the literature. In the newly obtained solution, the
expressions for the levels of consumption and capital stock are
equal to those found in the older solution, while the amount of
labor allocated to the production of physical capital and the level
of human capital are different from those values found in the older
solution. The existence of three closed form solutions is new to the
literature. We also show that these equilibria all converge to the
same long run balanced growth path.

\end{document}